\documentclass[10pt,twoside]{article}
\usepackage[a5paper,layoutwidth=148.5mm,layoutheight=7.77817in,
outer=10mm,inner=38pt,tmargin=10mm,bmargin=10mm,footskip=15pt,
]{geometry}
\usepackage[utf8]{inputenc}
\usepackage[T1]{fontenc}
\usepackage{lmodern}
\usepackage{fix-cm}
\usepackage{microtype}
\usepackage{textcomp}
\newcommand\unite{\text{\sffamily\large I}}

\newcommand\bruch{
}
\usepackage[comma]{natbib}

\setcitestyle{aysep=}
\usepackage[german,french]{babel}
\usepackage[french]{varioref}
\usepackage[autolanguage]{numprint}

\usepackage{graphicx}

\usepackage{amsmath}
\usepackage{amsfonts}

\usepackage{indentfirst}
\usepackage{framed}

\usepackage{titlesec}
\titleformat{\section}[runin]{\bfseries}{\textbf{\thesection.}}{\fontdimen2\font plus \fontdimen3\font minus \fontdimen4\font}{}
\titlespacing{\section}{0pt}{1.05ex plus.3ex minus.06ex}{2\fontdimen2\font plus 2\fontdimen3\font minus 2\fontdimen4\font}

\usepackage{xurl}
\PassOptionsToPackage{urlcolor=magenta,citecolor=blue,anchorcolor=blue,linkcolor=blue,menucolor=cyan}{hyperref}
\usepackage[pdftitle=,pdfauthor=Stefan\ Neuwirth,pdfdisplaydoctitle=true,pdflang=fr-FR,bookmarksopen=false,breaklinks=true,unicode=true,plainpages=false,hyperindex=true,pdfstartview=FitH,colorlinks=true,pdfpagelabels=true,bookmarksnumbered=true]{hyperref}
\ifx\texorpdfstring\undefined\newcommand\texorpdfstring[2]{#1}\fi
\ifx\pdfbookmark[3]{}\fi
\usepackage{doi}

\usepackage{framed}

\begin{document}

\title{\vspace{-10mm}«~C'est comme œuf~»\\{\large Paul Lorenzen et l'effondrement des preuves de non-contradiction}}
\author{Stefan Neuwirth}
\date{\vspace{-5mm}}
\maketitle

\section{}

Paul Lorenzen, mathématicien et philosophe du 20\textsuperscript{e} siècle, évoque le mois d'octobre 1947 comme la date d'une césure dans ses investigations tant mathématiques que philosophiques. J'ai découvert un écrit autographe daté du 15 octobre 1947 qui marque cette césure. J'en propose une traduction et un commentaire ligne à ligne, puis j'esquisse les circonstances de son écriture sur la base de la correspondance avec Paul Bernays, figure tutélaire de la logique mathématique et l'un de ses principaux interlocuteurs à cette époque, à l'instar de Heinrich Scholz et Oskar Becker. Une lettre de Lorenzen à Carl Friedrich Gethmann datée du 14 janvier 1988, dont je traduis un extrait, étoffe le récit de cette césure en montrant comment il s'imprègne des indications que lui fournissent ses correspondants pour les transmuter en une recherche absolument originale.

\section{}

Dans l'un des cartons déposés au Philosophisches Archiv de l'université de Constance par la fille de Lorenzen, noyé dans la correspondance, se trouve un feuillet jauni, isolé\footnote{Dans l'inventaire du Paul-Lorenzen-Nachlass, accessible sous \url{http://www.uni-konstanz.de/FuF/Philo/philarchiv}, il a le numéro PL~1-1-127 et est décrit (suite à l'intérêt que je lui ai porté) comme un «~feuillet recto-verso avec un texte programmatique, daté du 15.10.1947~», un mercredi.}. Il est manuscrit, d'une écriture soignée, porte un titre, un sous-titre, l'indication de l'auteur et du lieu, et se termine sur la date de la rédaction. Voici une traduction de cet écrit.

\begin{framed}\parindent0pt
  \hfil Pourquoi \textdied\qquad Vérité \textdied\medskip
  
  \hfil Une pièce ouvrée à la mémoire de Blaise Pascal.
  
  \hfil De Paul Lorenzen, Bonn.\bigskip

  Parler, c'est toujours faire une parabole.

  Question. Qu'est-ce que faire une parabole?

  Réponse. J'ai une montre. J'appelle celle-ci «~œuf~». Si j'examine le monde, je dis parfois: «~C'est comme œuf~». J'appelle ces mots combinés «~parabole~».

  Qui examine cette pièce ouvrée et la comprend, je dis de lui: «~il a vu l'œuf et la parabole~». Il peut palper les mots «~œuf~» «~parabole~» et il les a saisis.

  Les mots  «~œuf~» «~parabole~» «~palper~» «~saisir~» sont comme des pelotes. Les mots «~est~» «~et~» «~ne\dots\ pas~» sont comme de petits crochets, comme pour combiner des pelotes.

  Qui peut prolonger la pièce \unite, \unite\unite, \unite\unite\unite, \unite\unite\unite\unite\ a saisi le comptage.\bruch

  Je n'emploie pas «~pourquoi~» ni «~vérité~». J'abandonne ces mots aux barbares. J'emploie parfois «~d'où~», parfois «~comment~». Je distingue «~évident~» et «~persistant~». J'emploie souvent les mots:
\[
  \begin{tabular}{c}
    \begin{tabular}{cc}
      tête&cœur\\
      clair -- confus&bon -- mauvais
    \end{tabular}\\[\bigskipamount]
    vouloir -- craindre -- espérer -- haïr -- aimer\\[\smallskipamount]
    voir -- comprendre
  \end{tabular}
\]
  Depuis aujourd'hui, j'ai cessé d'employer le mot «~problème du continu~».

  J'aime pour mon langage le mot «~théorie~».\medskip

  \hfil 15/10/1947  
\end{framed}

\section{}

Je propose un commentaire ligne à ligne de cet écrit.
\begin{itemize}
\item Le titre contient deux mots, chacun suivi d'une croix latine, qui peut indiquer le décès. Déjà la présence de ce signe évoque le \emph{Mémorial} de Blaise Pascal et la nuit de feu du 23 novembre 1654. Le sous-titre assume et confirme cette référence. Voir le §~\ref{pascal}.
\item \emph{Pièce ouvrée [Werkstück]}. Ce mot, que je traduis littéralement, fait de cet écrit un bloc de granit, déjà dégrossi mais toujours à travailler. 
\item \emph{Parler [Sprechen]}. Par l'usage de l'infinitif, Lorenzen reporte l'attention du produit du langage sur l'acte de parole.
\item \emph{Montre}. \emph{Œuf}. L'œuf de Nuremberg [Nürnberger Ei] est parmi les premières montres à gousset jamais produites. Nommer la montre «~œuf~» corrobore une ressemblance visuelle et tactile entre la montre et l'œuf de poule.
\item \emph{Parabole [Gleichnis]}. C'est l'acte de parole «~C'est comme œuf [Das ist wie Ei]~» tout entier qui est une parabole. Le sujet «~ce~», pronom démonstratif, est au premier abord mystérieux: qu'est-ce que ce «~ce~»? De même, l'emploi du mot «~œuf~» sans article est au premier abord surprenant et déroutant: on pourrait y voir une référence à la matière de l'œuf, qu'on traduirait avec l'emploi du partitif: «~de l'œuf~». J'interprète au contraire cette phrase ainsi:
  \begin{itemize}
  \item «~Ce~» réfère à l'acte d'examiner le monde, c'est-à-dire d'y discerner des objets, de leur accorder une individualité par des noms qui dans cet examen jouent le rôle de noms propres.
  \item «~Œuf~» réfère à l'acte d'appeler la montre «~œuf~» exactement comme le ferait un nom propre. Cet acte individué est ici la parabole de l'acte générique d'examiner le monde.
  \end{itemize}
\item \emph{Examiner [besehen]}. \emph{Comprendre [verstehen]}. \emph{Voir [ersehen]}. C'est en « examinant » et en « comprenant » qu'on « voit » au sens  de « connaitre ».
\item \emph{Pelotes [Bälle]}. On pourrait traduire ici par «~balles~», mais j'ai préféré «~pelotes~» pour renforcer l'image concrète du crochet qui accroche les pelotes entre elles.
\item «~Est~» «~et~» «~ne\dots\ pas~» sont des termes d'un langage élémentaire (Lorenzen dira «~protolangage~») élaboré afin d'organiser l'examen du monde.\bruch
\item \unite, \unite\unite, \unite\unite\unite, \unite\unite\unite\unite. Je vois \unite, \unite\unite, \unite\unite\unite, \unite\unite\unite\unite\ comme la parabole de la saisie que l'acte qui mène de \unite\ à \unite\unite\ ressemble à l'acte qui mène de \unite\unite\ à \unite\unite\unite\ et de \unite\unite\unite\ à \unite\unite\unite\unite; prolonger la pièce [das Stück verlängern] est pour moi la parabole de l'acte de répéter. Voir le §~\ref{unite}.
\item \emph{Pourquoi}. \emph{Vérité}. Selon moi, ces mots réfèrent à la crise que Lorenzen traverse à ce moment précis et que je décris plus loin. La question «~pourquoi?~» engage une régression à l'infini puisqu'elle ouvre la voie à la question du pourquoi du pourquoi, etc. La «~vérité~» réfère à un deuxième monde, le lieu des vérités, alors qu'il y a un seul monde à examiner.
\item \emph{Barbares}. Le mot est très fort et en 1947, il évoque irrésistiblement la Seconde Guerre mondiale et les horreurs du nazisme. Il est difficile de corroborer cette lecture. À minima, l'emploi des mots «~pourquoi~» et «~vérité~» a mené Lorenzen dans une impasse dont il décide de s'échapper.
\item \emph{D'où}. Ce mot réfère à une causalité concrète dans l'espace et le temps, relative à un avant et un après.
\item \emph{Comment}. Ce mot réfère à une description précise des actes matériels et de leur rapport avec les actes de parole.
\item \emph{Évident [ersichtlich]}. \emph{Persistant [beständig]}. C'est par une certaine stabilité et régularité que le monde s'offre à notre examen, et on constate cette régularité dans notre connaissance matérielle; Lorenzen distingue ceci de l'évidence propre à la connaissance formelle. Voir le §~\ref{ueber}.
\item \emph{clair -- confus [klar -- trüb]}, etc. Tous les mots de cet arrangement réfèrent à des expériences vécues et renvoient ainsi à la philosophie de la vie de Wilhelm Dilthey; voir le §~\ref{leben}. L'opposition de la tête au cœur fait penser à Pascal et à son fameux aphorisme «~Le cœur a ses raisons que la raison ne connait point~». Mais elle renvoie plus directement à la citation que Heinrich \citet{scholz45} extrait du fragment « Grandeur » des \emph{Pensées} \citep{pascala} pour établir l'« \emph{intuitionnisme} pascalien, dans la langue de Pascal la logique du cœur »: «~Nous connaissons la vérité non seulement par la raison mais encore par le cœur, c’est de cette dernière sorte que nous connaissons les premiers principes [\dots{}] Et c’est sur ces connaissances du cœur et de l’instinct qu’il faut que la raison s’appuie et qu’elle y fonde tout son discours [\dots{}] Les principes se sentent, les propositions se concluent.~»
\item \emph{Problème du continu [Kontinuumproblem]}. Ce problème est le premier des 23 fameux problèmes que David Hilbert propose au congrès international des mathématiciens à Paris en 1900. Lorenzen a appris quelques mois auparavant, comme il l'évoque dans sa lettre à Bernays du 8 avril 1947 (ETH-Bibliothek, Hochschularchiv, Hs~975\string:\allowbreak2952), que Kurt Gödel a conçu un «~modèle~» formé des «~ensembles constructibles~» qui satisfait tous les axiomes de la théorie des ensembles infinis, n'est pas plus contradictoire que celle-ci et répond positivement au problème du continu. Par sa décision de cesser d'employer ce mot, il récuse l'affirmation qu'il s'agirait d'un problème bien posé et s'en détourne.
\item \emph{Théorie}. L'écrit se termine par l'appréciation positive de ce mot, peut-être parce que la théorie est le lieu de la constatation et de l'étude des règles. Ce mot est pourtant rare dans les textes de Lorenzen de cette époque: il préfère travailler sur les formalisations de théories, qu'il appelle «~calculs~».\bruch
\end{itemize}

\section{}

Qui est Paul Lorenzen en 1947? Il est au premier chef mathématicien. Il n'a cessé de faire des mathématiques tout au long de la Seconde Guerre mondiale, pendant laquelle il a d'abord été simple soldat, puis s'est montré incapable d'obtempérer silencieusement dans l'emploi que son directeur de thèse Helmut Hasse lui avait trouvé dans la marine, en affichant un «~comportement non militaire~» (lettre de Käthe Lorenzen à Hasse datée du 7 mai 1941, \foreignlanguage{german}{Helmut-Hasse-Nachlass, Niedersächsische Staats- und Universitätsbibliothek Göttingen}, Cod.\ Ms.\ H. Hasse 1\string:\allowbreak1021). Il s'est finalement retrouvé enseignant à l'école de marine de Wesermünde. Au retour de la guerre, il reprend le poste d'assistant auprès du professeur Wolfgang Krull à Bonn auquel il a été nommé en 1939. Il est habilité en 1946 avec un mémoire sur la théorie des idéaux pour les groupes préordonnés. Mais il est en rupture avec Krull, suite aux obstacles que ce dernier n'a cessé d'opposer à cette habilitation pourtant nécessaire pour envisager une carrière scientifique parce qu'«~une nouvelle leçon ne lui ferait pas de mal~» (carte envoyée par Krull à Hasse datée du 6 février 1944, éditée dans \citealt[\S~1.89]{roquette04}). Je conjecture qu'ils ont voulu lui infliger cette leçon pour qu'il se rende compte de quel genre d'homme le Troisième Reich avait besoin et parce qu'ils craignaient de commettre une erreur politique s'ils lui ouvraient la voie pour une carrière universitaire. La correspondance de Lorenzen n'en fait pas montre, mais celui-ci est d'un caractère ombrageux, entier; il est certain de son chemin, radical dans sa pensée, réfractaire par nature; Hasse et Krull le jugent «~fier~» et «~présomptueux~».

\section{}
\label{leben}

Lorenzen n'est pas seulement mathématicien: il a fait un cursus parallèle en philosophie pendant ses années d'études à Göttingen, auprès de Josef König et de Fritz Gebhardt qui se placent dans le courant de la philosophie de la vie de Wilhelm Dilthey. \citet{MR0284309} écrit que «~La réflexion sur l'herméneutique, c.-à-d.\ sur la théorie de la compréhension des actes humains, en particulier d'assertions parlées et écrites, a conduit de Schleiermacher à Dilthey chez celui-ci à l'assertion remarquable suivante: ``La connaissance ne peut pas reculer en deçà de la vie\footnote{Lorenzen ne donne pas la source de cette assertion; dans «~Zur Weltanschauungslehre~» (publié en 1931 dans le tome~VIII des \emph{Gesammelte Schriften}, page~180) on trouve «~Hinter das Leben kann das Erkennen nicht zurückgehen~».}''~».

\section{}

Voici comment Lorenzen raconte lui-même son cheminement dans une lettre à Carl Friedrich \citet{gethmann91}.

\begin{framed}
  Mais l'évolution vers l'approche constructive de critique du langage s'accomplit hors de la discipline «~philosophie~» par mes travaux algébriques qui s'employaient à un problème qui avait \emph{formellement} la même structure que le problème de la non-contradiction du calcul logique classique (s'il est utilisé pour des domaines infinis comme dans l'arithmétique). J'ai suivi à Göttingen, vers 1938, un exposé de Gentzen sur sa preuve de non-contradiction. Mes travaux algébriques traitaient de la «~théorie des idéaux~» dans les corps de nombres algébriques. Il s'avère que dans ces corps, une relation de divisibilité ($\alpha$ divise~$\beta$) est définie sans que (comme pour les nombres rationnels) un «~plus grand commun diviseur~» (p.g.c.d.)\ et un «~plus petit commun multiple~» (p.p.c.m.)\ doive toujours exister. La théorie des idéaux \emph{construit} de tels p.g.c.d.\ et p.p.c.m.\ comme des éléments idéaux.

  Dans la logique, une relation d'implication est définie pour les propositions élémentaires (par des règles de prédicateurs, comme on dira plus tard). \bruch Elle est réflexive et transitive, comme la relation de divisibilité. P.g.c.d.\ et p.p.c.m.\ correspondent formellement à des conjonctions et adjonctions (pour une infinité de propositions élémentaires -- comme elles résultent de formules propositionnelles par substitution~--, ce sont les quanteurs). Cela prit jusqu'en 1947 (entre temps il y avait toujours la guerre -- et quelques visites auprès de Scholz à Münster) jusqu'à ce qu'il devint clair pour moi que des preuves de non-contradiction ne pouvaient pas être menées dans des «~métalangages~» (parce que la logique du métalangage requiert alors à nouveau une preuve de non-contradiction, et ainsi à l'infini), que tout au contraire les propositions logiques complexes sont \emph{construites} par nous.

  L'effondrement des preuves de non-contradiction et des métalangages et le recommencement avec le protolangage, cela eut lieu en octobre 1947. Seul O. Becker fut compréhensif à cet égard; il connaissait le problème par Husserl. Je n'ai cependant jamais compris pourquoi ni Husserl ni Heidegger (ni Becker lui-même) n'ont sérieusement \emph{construit} des langages scientifiques. Becker disait alors toujours: s'il vous plait, concédez au vieux sceptique que je suis de ne pas partager vos espoirs qu'un «~ortholangage~» pourrait être imposable.
\end{framed}

\section{}

Comme il le raconte, Lorenzen poursuit parallèlement à ses recherches algébriques des recherches logiques, bien qu'il ait été averti par Krull que ces travaux ne vaudraient pas pour son habilitation (lettre du 19 février 1944, PL~1-1-143, éditée dans \citealt[p.~229-230]{neuwirthkonstanz}). Le besoin de sonder les fondations de son discours algébrique l'amène à réfléchir sur ce que c'est qu'être un idéal.
\begin{itemize}
\item Pour Eduard Kummer, les nombres idéaux sont de simples façons de parler auxquelles on n'attribue pas d'existence hors du discours: ce qui existe, c'est la relation de divisibilité entre nombres idéaux. En d'autres mots, Kummer introduit les nombres idéaux comme Euclide introduit les rapports de grandeurs dans le livre~V des \emph{Éléments}: ceux-ci ne sont pas définis; c'est la relation de comparaison entre rapports de grandeurs qui l'est.
\item Pour Richard Dedekind, il est nécessaire d'attribuer une existence absolue et réelle aux rapports de grandeurs et aux nombres idéaux. Il y parvient en se servant de manière cruciale de la théorie des ensembles infinis, qui postule à l'avance l'existence de domaines infinis (en particulier de l'ensemble des sous-ensembles de l'ensemble~$\mathbb Q$ des nombres rationnels ou d'un corps de nombres algébriques) pour permettre à un raisonnement mathématique d'y faire exister l'extension des concepts considérés (c'est-à-dire les sous-ensembles de tous les éléments qui tombent sous le concept). Cette postulation ne répond en rien à la question de la construction effective de ces objets.
\item Lorenzen propose une autre voie: il ne conçoit pas les idéaux comme certains ensembles infinis dont l'existence est simplement postulée, mais comme des objets formels construits par des règles de formation et des règles de raisonnement qui résultent immédiatement de leur signification, comme suit. Étant donnés~$a$ et~$b$ pris dans un corps de nombres algébriques, on n'y trouve pas toujours un~$d$ qui serait le plus grand parmi tous les nombres~$x$ qui divisent à la fois~$a$ et~$b$, c'est-à-dire tel que
  \[\text{si $x$ divise $a$ et $x$ divise $b$, alors $x$ divise~$d$.}\]
  Lorenzen propose de construire~$d$ comme un objet nouveau, la combinaison de trois signes $a\land b$, introduite avec la règle de raisonnement\bruch
  \[\text{si $x$ divise $a$ et $x$ divise $b$, alors $x$ divise~$a\land b$,}\]
  qu'on rajoute aux règles de calcul dans le corps de nombres. Ceci \emph{construit} le p.g.c.d.\ de~$a$ et~$b$ si on peut prouver qu'on a bien
  \[\text{si $x$ divise~$a\land b$, alors $x$ divise $a$ et $x$ divise $b$,}\]
  c'est-à-dire qu'on dispose de la règle de raisonnement qui élimine la combinaison~$a\land b$. Pour l'établir, Lorenzen montre qu'elle est «~admissible~», c'est-à-dire que toute preuve qui aboutit à son hypothèse «~$x$ divise~$a\land b$~» peut être transformée en une preuve qui aboutit à sa conclusion «~$x$ divise~$a$ et $x$ divise~$b$~».
\end{itemize}

Ce faisant, Lorenzen se rend compte que cette question a son pendant exact dans l'élaboration du calcul logique, où il s'agit de répondre à la question: qu'est-ce que c'est qu'être une conjonction, une quantification? Cette prise de conscience est documentée dans sa correspondance avec Krull, en particulier dans sa lettre du 13 mars 1944: «~À la suite d'une investigation algébrique des demi-treillis orthocomplémentés, j'essaie maintenant d'élucider le rapport de ces questions avec la non-contradiction de la logique classique. Je veux faire cet essai même au risque que vous ne teniez pas compte du résultat éventuel -- parce que je ne peux pas faire autrement que de considérer de telles questions comme des questions d'un intérêt mathématique primordial~» (PL-1-1-131, éditée dans \citealt[p.~230-232]{neuwirthkonstanz}).

\section{}
\label{unite}

Cependant, son investigation le pousse plus loin: alors que sa preuve de non-contradiction de 1944 ramène la non-contradiction de la théorie élémentaire des nombres à la non-contradiction du concept même de nombre entier (voir \citealt{coquandneuwirth19}), c'est l'élaboration de ce concept qu'il veut comprendre.

Son intention est de se servir de la technique qu'il a développée au sujet de la logique pour la poursuivre au sujet du nombre. Il tient Bernays au courant de ses réflexions en lui envoyant régulièrement des nouvelles:
\begin{itemize}
\item «~Je me suis entre temps efforcé de formaliser les déductions utilisées dans la preuve de non-contradiction.
Il s'agit d'une restriction du calcul intuitionniste par laquelle les preuves indirectes sont exclues. La sémiotique est édifiée avec cette logique\footnote{«~Ich habe mich inzwischen darum bemüht, die im Wf.beweis benutzten
Schlüsse zu formalisieren. Es handelt sich dabei um eine Einschränkung
des intuitionistischen Kalküls, durch die indirekte Beweise ausgeschlossen werden. Mit dieser Logik wird die Semiotik aufgebaut~» (lettre du 21 février 1947, Hs~975\string:\allowbreak2950).}.~»
\item Il s'exprime de manière semblable dans sa lettre du 4 mai 1947 (Hs~975\string:\allowbreak2953).
\item «~Plus fondamentale me parait la possibilité d'éviter par la méthode des métacalculs toute référence au ``contenu'', à l'``intuitif'', toute ``évidence finitaire'' (que je veux qualifier d'``actualisme finitaire\footnote{«~\foreignlanguage{german}{Wesentlicher scheint
mir die Möglichkeit, durch die Methode der Metakalküle alles
``Inhaltliche'', ``Anschauliche'', alle ``finite Evidenz'' (was ich als
``finiten Aktualismus'' bezeichnen möchte) zu vermeiden}~» (lettre du 27 mai 1947, PL~1-1-106).}'').~»
\end{itemize}\bruch

Les tentatives entreprises deviennent de plus en plus techniques et ont trait à la «~logique de conséquence~». Parallèlement, sa compréhension du nombre entier s'approfondit. En particulier, la liste \unite, \unite\unite, \unite\unite\unite, \unite\unite\unite\unite\ de l'écrit fait écho à une réflexion annoncée dans la lettre du 4 mai 1947 à Bernays: «~L'arithmétique élémentaire aussi, je la fonderais par construction et non par évidence intuitive\footnote{«~\foreignlanguage{german}{Auch die elementare Arithmetik würde ich durch Konstruktion begründen, nicht durch anschauliche Evidenz}~»). La réflexion est poursuivie dans le manuscrit joint à cette lettre, «~\foreignlanguage{german}{Konstruktive Begründung der Arithmetik [Fondement constructif de l'arithmétique}]~».}.~»

Il construit le nombre entier comme suit: le premier, appelé l'\emph{unité}, est le signe~\unite; étant donné un nombre déjà construit~$x$, on construit son successeur comme la combinaison de~$x$ avec le signe~\unite, notée~$x\unite$. Si on prend pour~$x$ l'unité, on construit la combinaison~\unite\unite, appelée le \emph{deux}; si on prend pour~$x$ le deux, on construit la combinaison~\unite\unite\unite, appelée le \emph{trois}; et ainsi de suite.

\section{}

Cependant, Bernays lui rappelle dans sa lettre du 1\textsuperscript{er} septembre 1947 que la constatation de la non-contradiction repose sur certains «~discernements~»: «~Il me parait très douteux qu'il soit possible de se passer de l'évidence intuitive dans la théorie de la démonstration. Si on exécute des constructions, cela demande encore certains discernements, pour qu'on puisse tirer de ceux-ci une conséquence d'un caractère général comme l'est en particulier toute constatation de non-contradiction d'une théorie formelle\footnote{«~Ob es möglich ist, in der Beweistheorie ohne anschauliche Evidenz auszukommen, erscheint mir als sehr fraglich. Wenn man Konstruktionen ausführt, so bedarf es doch immer noch gewisser Einsichten, um aus diesen ein Ergebnis von allgemeinem Charakter zu gewinnen, wie es doch insbesondere jede Feststellung von Widerspruchsfreiheit einer formalen Theorie ist~» (PL~1-1-112).}.~»

\section{}
\label{ueber}

C'est peut-être cette réserve de Bernays qui pousse Lorenzen à mener son investigation plus loin et à réfléchir sur la nature de ces «~discernements~», jusqu'à aboutir à l'écrit présenté ici. Il se met alors à la rédaction d'un manuscrit intitulé «~Sur la connaissance [Über das Erkennen]~» qu'il joint à sa lettre du 1\textsuperscript{er} novembre 1947 (Hs~975\string:\allowbreak2956), dans laquelle il raconte: «~La santé de ma femme est meilleure ces temps-ci qu'elle ne l'a été depuis longtemps -- c'est pourquoi c'est moi qui suis soigné au moyen de votre colis [Bernays envoie régulièrement des colis alimentaires à la famille Lorenzen ainsi qu'à d'autres contacts], car je dois avouer que l'élaboration du manuscrit ci-joint m'a maintenu plus que sous tension 15 jours durant -- pendant lesquels j'ai dormi quelque peu tout au plus trois nuits\footnote{«~Meiner Frau geht es zur Zeit gesundheitlich besser als seit langem --~daher werde ich zur Zeit mit den Mitteln Ihres Paketes gepflegt, denn ich muß gestehen, daß die Erarbeitung des beiliegenden Manuskripts mich 14~Tage lang --~während der ich höchstens drei Nächte etwas geschlafen habe~-- mehr als in Spannung gehalten hat.~»}.~»

Voici les premières lignes de ce manuscrit.
\begin{framed}
  Notre connaissance est une unité de matière et de forme, de connaissance matérielle et formelle.\bruch

  Le langage de la connaissance matérielle consiste en des noms. Nous distinguons les noms propres comme Pascal, Goethe, Einstein\dots\ des noms communs comme maison, rouge, expérience vécue, bien.

  Nous connaissons des ressemblances et des différences\footnotemark.
\end{framed}
\footnotetext{«~Unsere Erkenntnis ist eine Einheit von Stoff und Form, von materialer
und formaler Erkenntnis.

«~Die Sprache der materialen Erkenntnis besteht aus Namen. Wir unterscheiden Einzelnamen wie Pascal, Goethe, Einstein, \dots\ von den Gemeinnamen wie Haus, rot, Erleben, gut.

«~Wir erkennen Ähnlichkeiten und Unterschiede.~»}

Il se rend compte que l'usage même de la parole est soumis à des règles et qu'en dernier lieu il repose sur des jugements de ressemblance et de dissemblance. Il y revient dans sa lettre du 10 décembre 1947, dans laquelle il admet que le manuscrit ne restitue pas sa pensée: «~De ce qui est écrit dans les feuillets ``Sur la connaissance'', on ne peut pas extraire ce que je voulais dire même avec beaucoup de bonne volonté\footnote{«~dass man aus dem, was in den Blättern ``Über das Erkennen'' steht, selbst bei gutem Willen nicht das herauslesen kann, was ich meinte~» (Hs~975\string:\allowbreak2957).}.~» Il y évoque la «~critique de la connaissance [\dots] qui ne laisse subsister que ressemblance et différence (mais ce ``ne\dots\ que'' est beaucoup -- il ne doit exprimer aucun jugement de valeur\footnote{«~Erkenntniskritik [\dots], die material „nur“ Ähnlichkeit und Unterschied bestehen lässt
(dies „nur“ ist aber sehr viel – es soll keinerlei Wertung ausdrücken)~».})~».

Un but affirmé de sa démarche est d'élaborer un «~langage contrôlable~» (expression qu'il attribue à Heinrich Scholz) pour «~que dans les discussions sur les mathématiques on doive mettre le doigt sur l'endroit où quelque chose devrait être faux\footnote{«~Das Ziel meiner Untersuchungen ist der Zustand, den Hilbert vorausgesehen hat: ``daß man bei allen Diskussionen über die Mathematik den Finger auf die Stelle legen muß, wo etwas falsch sein sollte'' (ungenau zitiert). Es ist dasselbe Ziel, das Scholz als ``kontrollierbare Sprache'' bezeichnet~» (lettre du 8 février 1948, Hs~975\string:\allowbreak2959). Lorenzen dit citer Hilbert de manière imprécise; \citet[p.~158]{hilbert22} avance qu'«~il devient possible de contrôler en quelque sorte directement chaque formule de la démonstration, c'est-à-dire qu'on est en mesure d'établir si elle est ``vraie'' ou ``fausse'' en un sens déterminé qu'il resterait à décrire avec précision~»; mais Lorenzen se réfère sans doute à \citealt[p.~140-141]{hilbert29}: «~La démonstration de la consistance permet en même temps, chaque fois qu'on est en présence d'une démonstration conduisant à un résultat faux, de trouver l'endroit de l'erreur~»; «~Quiconque prétend me réfuter, doit, selon un usage qui en mathématiques a toujours prévalu dans le passé et continuera de prévaloir à l'avenir, me montrer exactement où je suis censé m'être trompé~».}~». Sa lettre du 1\textsuperscript{er} novembre corrobore le récit qu'il fait à Gethmann: «~Je n'ai guère pu trouver de compréhension pour mon langage contrôlable ici à Bonn\footnote{«~Für meine kontrollierbare Sprache habe ich hier in Bonn kaum Verständnis finden können.~»}.~»

\section{}
\label{pascal}

Le premier texte que Lorenzen écrit et publie après la césure du 15 octobre 1947 est son \emph{Einführung in die Logik [Introduction à la logique]} dont il joint un exemplaire à sa lettre à Bernays du 8 février 1948 et dont il publie une recension \citep{MR0033260}. Il y écrit: «~[\dots]\ une règle n'est introduite que si sa signification est clairement manifeste. [\dots] toute règle en question est seulement une règle d'action, c.-à-d.\ une règle selon laquelle on doit agir avec des signes (figures). Si on pose la question de la signification d'une telle règle, on pose la question de la signification et de la finalité de l'action selon la règle. [\dots] La question de la signification de la construction de ces domaines [des nombres et des égalités arithmétiques] est résolue par l'indication de leur utilité dans la vie quotidienne.~»

L'article de \citeauthor{MR0284309} de \citeyear{MR0284309} revient plus spécifiquement sur le rapport qu'il établit entre sa pensée et celle de Pascal, en particulier au sujet de la régression à l'infini, en citant le passage suivant du fragment «~Disproportion de l'homme~» dans les \emph{Pensées} \citep{pascalb}: «~C'est ainsi que nous \bruch voyons que toutes les sciences sont infinies en l'étendue de leurs recherches, car qui doute que la géométrie par exemple a une infinité d'infinités de propositions à exposer. Elles sont aussi infinies dans la multitude et la délicatesse de leurs principes, car qui ne voit que ceux qu'on propose pour les derniers ne se soutiennent pas d'eux-mêmes et qu'ils sont appuyés sur d'autres qui en ayant d'autres pour appui ne souffrent jamais de dernier.~» Notons que Heinrich Scholz aussi a étudié l'œuvre de Pascal; peut-être est-ce d'ailleurs lui qui l'a fait découvrir à Lorenzen.

\section{}

L'écrit présenté ici, dont les successeurs de Paul Lorenzen comme Kuno Lorenz et Christian Thiel ignoraient jusqu'à l'existence, n'a peut-être jamais quitté les mains de son auteur. L'a-t-il seulement montré à Oskar Becker, Fritz Gebhardt, Josef König ou Heinrich Scholz?

Il témoigne d'un événement qu'il qualifie d'effondrement: Lorenzen est amené à abandonner une posture qui s'est révélée erronée et une entreprise qui s'est avérée être une impasse. Le titre porte les marques de l'abandon de deux exigences:
\begin{itemize}
\item celle de la recherche de la vérité comme d'un but en soi, comme si elle existait en soi et qu'il s'agissait  de s'en approcher;
\item celle de la recherche des causes, comme si elle pouvait aboutir une fois pour toutes à un fondement de notre connaissance.
\end{itemize}

Il évoque trois problèmes mathématiques:
\begin{itemize}
\item celui de la signification des connecteurs logiques comme
  les mots «~est~» «~et~» «~ne\dots\ pas~»;
\item celui de la signification du comptage qui donne lieu à \unite, \unite\unite, \unite\unite\unite, \unite\unite\unite\unite;
\item celui de la signification du problème du continu.
\end{itemize}
Les deux premiers sont reconnus comme des problèmes en lien avec la vie et avec notre volonté d'organiser la connaissance; le troisième est rejeté comme un problème inerte, dont la signification ne se fonde pas dans la vie connaissante, mais dans l'expérience de laboratoire d'un jeu axiomatique (la théorie des ensembles) dont les règles ont été établies en opérant subrepticement une rupture irréparable avec la vie. La réparation que Hilbert espérait était une preuve de non-contradiction des règles de ce jeu, or ce moyen est justement celui dont Lorenzen constate l'effondrement.

\section{}

Lorenzen accorde une place centrale au langage et décrit la connaissance comme l'articulation de nos facultés de nommer des choses et de combiner des noms avec des expériences vécues qui fondent la signification de ces combinaisons. Il apparait par cet écrit comme un Don Quichotte qui a entrepris de poursuivre le programme de Hilbert, fonder les mathématiques par une preuve de non-contradiction, jusqu'au cœur de ce programme, qui a fait l'expérience à son corps défendant de la régression à l'infini, et qui a dû se résoudre à refonder ce programme sur une articulation entre la théorie et la vie. 

Cet écrit a pour point de départ le moment poétique de la connaissance. Son auteur s'y révèle poète. Sa difficulté patente à trouver les mots, à extérioriser le cheminement de sa pensée, à partager \bruch ses réflexions, son gout (dénoncé dès la rédaction de sa thèse par Hasse et Krull, voir les lettres qu'ils s'échangent les 21 et 31 mai 1938, \citealp[§~1.35, 1.36]{roquette04}) pour une écriture dense et ramassée obligent le lecteur à parcourir à son tour ce cheminement en n'y étant pas vraiment guidé, mais en étant accueilli à l'arrivée par des indications indubitables qu'il a bien gravi le sommet que l'auteur avait gravi auparavant: son expression elliptique est en effet d'une précision et d'une maitrise telles qu'alors la confusion du lecteur se dissipe et laisse place à la clarté, comme j'en ai fait l'expérience avec l'affirmation «~C'est comme œuf~».

\section{}

Terminons sur un extrait de la lettre de Lorenzen à Scholz datée du 4 mai 1947: «~Je lis en ce moment à nouveau dans la ``Théorie des idées de Platon'' de Natorp. Par suite, il me semble que le constructivisme ne contredit pas du tout la conception de Platon. Car, que les lois soient quelque chose d'étant, que d'elles ne peut être connu que ce qu'elles sont en elles-mêmes, qu'elles seules sont donc véritablement étantes, que nous ne pouvons connaitre qu'en agissant selon des lois, que par notre participation aux lois -- pourquoi devrais-je le contester. C'est plutôt la description la plus exacte de ma conception\footnote{«~\foreignlanguage{german}{Ich lese zurzeit wieder in ``Platos Ideenlehre'' von Natorp. Danach scheint mir der Konstruktivismus durchaus nicht der Platonischen Auffassung zu widersprechen. Denn, dass die Gesetze etwas
Seiendes sind, dass nur von ihnen erkannt werden kann, was sie an
sich selbst sind, dass sie also allein wahrhaft seiend sind, dass
wir nur durch gesetzmässiges Handeln, durch unser Teilhaben an den
Gesetzen erkennen können --~warum sollte ich dem widersprechen. Das ist vielmehr die genaueste Beschreibung meiner Auffassung}~» (PL~1-1-119).}~».
On peut lire l'extrait suivant de la préface de Scholz aux \emph{Grundzüge der
mathematischen Logik} (\citeyear{scholzhasenjaeger61}) comme une réponse à Lorenzen: «~Dans
le [cas d'une logique stipulée ontologiquement], les présuppositions
ontologiques sont tellement essentielles qu'on peut se passer par
principe d'actions au sens d'opérations. Il est essentiel de le savoir. La
logique présentée ici est indépendante de l'existence d'êtres capables
d'exécuter une quelconque action prescrite. C'est seulement dans cette
absence de présupposés que son caractère platonicien est entièrement
mis en valeur. Dans le [cas d'une logique stipulée constructivement]
c'est l'inverse: ici les opérations sont au fondement. La question de la
mesure dans laquelle une ontologie peut être abstraite d'elles à postériori
est secondaire. De toute façon on n'a pas besoin d'ontologie pour se
constituer. C'est le discernement présupposé des actions requises qui
est décisif.~»\medskip

\emph{Je remercie Agnès Leval et Philippe Séguin pour leur relecture. Je remercie Laurence Dahan-Gaida pour son invitation à la journée « Eurêka ! » qui a suscité cette recherche.}\clearpage

\noindent\includegraphics[width=\textwidth]{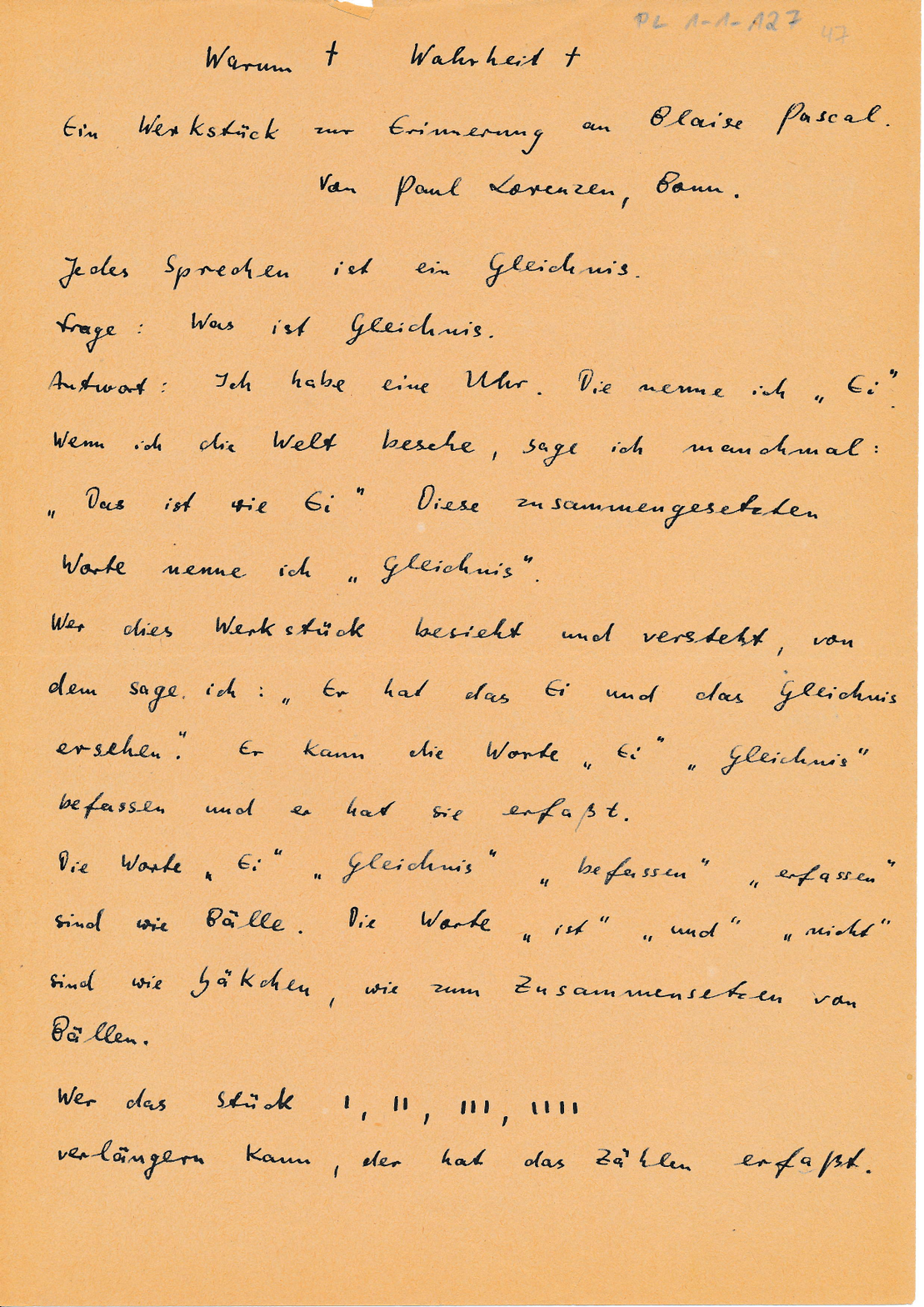}

\noindent\includegraphics[width=\textwidth]{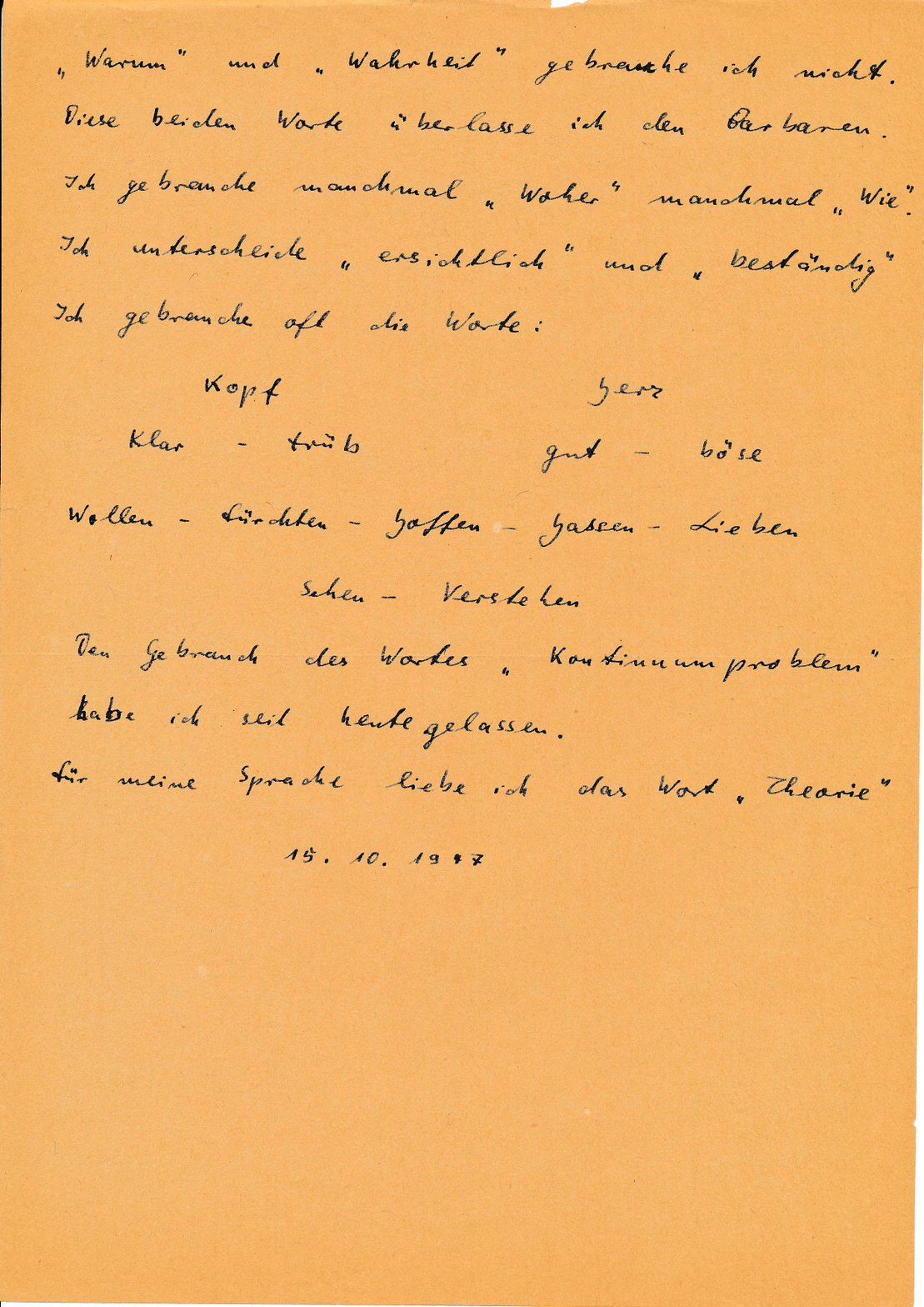}


\end{document}